\documentclass[a4paper] {article} [11pt]
\usepackage{amsmath,amsfonts,amssymb,amsthm,amscd,nccmath,verbatim}
\usepackage{graphicx}
\usepackage[numbers]{natbib}
\usepackage{color}
\usepackage{yfonts}
\usepackage[english]{babel}

\newtheorem{Theorem}{Theorem}[section]
\newtheorem{thm}[Theorem]{Theorem}

\newtheorem{prop}[Theorem]{Proposition}
 
\newtheorem{defn}[Theorem]{Definition} 

\newtheorem{rmk}[Theorem]{Remark}
\newtheorem{rem}[Theorem]{Remark}

\newtheorem{lem}[Theorem]{Lemma}

\newtheorem{cor}[Theorem]{Corollary}
\newtheorem{Fact}[Theorem]{Fact}

\newtheorem{ex}[Theorem]{Example}

\newtheorem{qu}[Theorem]{Question}

\newtheorem{example}[Theorem]{Example}

\newcommand{\QED}{\hfill \raisebox{-1em}{$\Box$}}
\newcommand{\mbb}{\mathbb}
\newcommand{\mc}{\mathcal}
\newcommand{\ra}{\rightarrow}
\renewcommand{\iff}{\leftrightarrow}

\title{The model theory of commutative near-vector spaces}
\author{Karin-Therese Howell, Charlotte Kestner}
\date{}
\begin{document}
\maketitle

\begin{abstract}
In this paper we study  near-vector spaces over a commutative $F$ from a model theoretic point of view. In this context we show regular near-vector spaces are in fact vector spaces. We find that near-vector spaces are not first-order axiomatisable, but that finite block near-vector spaces are. In the latter case we establish quantifier elimination, and that the theory is ``controlled'' by which elements of the pointwise additive closure of $F$ are automorphisms of the near-vector space.
\end{abstract}

\section{Introduction}

In \cite{Andre}, the concept of a vector space, is generalised by Andr\'e to a structure comprising a bit more non-linearity, the so-called near-vector space. This consists of a commutative group $V$ with a set $F$ of endomorphisms of $V$. In \cite{vdWalt2} Van der Walt showed how to construct an arbitrary finite-dimensional near-vector space, using a finite number of near-fields, all having isomorphic `multiplicative' semigroups. In \cite{HowM1} this construction is used to characterize all finite-dimensional near-vector spaces over $\mathbb F_p$, where $p$ is a prime. These results were extended to all finite-dimensional near-vector spaces over arbitrary finite fields in \cite{HowM2}.  Near-vector spaces have been used in many applications,
including in cryptography \cite{NT_NVS} and in interesting examples of classes of planar near-rings \cite{planar}.
Our aim in this paper is to use model theory to add to the theory and understanding of near-vector spaces.

In section 2 the basic facts and definitions are established. We attempt to give sufficient background for the paper to be self contained for both those near-vector space experts that are new to model theory and vice versa.
In section 3 we consider commutative near-vector spaces (i.e. composition of functions in $F$ is commutative). Although there are many examples where this is not the case
, there are several nice conclusions that can be drawn in the case where $F$ is commutative. We show that regular near-vector spaces are in fact vector spaces over $(F, \circ, +_u),$ where here the addition is that induced by any member of the quasi-kernel. This is used to show that any commutative near-vector space decomposes into `blocks', each of which will be a vector space over a field whose base set is $F$. We give a surprising example, where the induced fields have different characteristic and also clarify the statement in \cite{Andre}, regarding when a near-vector space is a vector space.




Section 4 examines the near-vector space as a one sorted structure (in the language of modules), and establishes some basic model theoretic properties of these structures. We show that near-vector spaces are not first-order axiomatisable in this language, but that those with finitely many `blocks' are.
We establish quantifier elimination in the latter case, under the assumption that $F$ is commutative. We define $\bar{F}$ the closure of $F$ under pointwise addition, and note that if $F$ is commutative then a near-vector space is in fact a module over $\bar{F}$, showing that these structures must in fact be stable (in the model theoretic sense). We show that for a `finite block' near-vector space the theory of $V$ is `controlled' by the set $\bar F \cap Aut(V)$
and use this to show that `finite block' near-vector spaces are totally transcendental with Morley rank equal to the number of blocks.

The authors would like to thank Gareth Boxall for some extremely valuable comments and contributions. We would also like to thank Lorna Gregory for some very helpful conversations around the model theory of modules.

\section{Preliminaries}

This paper establishes the basic model theoretic facts of near-vector spaces. We give both an introduction to near-vector spaces and one to model theory. The aim of this introduction is to provide the necessary background to make the paper accessible to a very general audience. 

\subsection{Near-vector spaces}

\begin{defn} (Definition 1.1., \cite{Andre})
An \textit{F-group}  is a structure $(V, F)$ which satisfies the following conditions:
\newline ($F_{1}$) $(V, +)$ is a group and $F$ is a set of endomorphisms of $V$;
\newline ($F_{2}$) The endomorphisms $0$, $1$ and $-1$, defined by $0x = 0$, $1x = x$ and $(-1)x = -x$ for each $x \in V$, are elements of $F$;
\newline ($F_{3}$) $F^{*} := F \backslash \{0\}$ is a subgroup of the group of automorphisms of ($V$, $+$);
\newline ($F_{4}$) If $\alpha x = \beta x $ with $x \in V$ and $\alpha, \beta \in F$, then $\alpha = \beta$ or $x = 0$, i.e. $\!F$ acts \textit{fixed point free} (fpf) on $V$.
\end{defn}

\medskip
\begin{rmk}
(a) $(V, +)$ is abelian, since by ($F_{2}$): 
\newline  $\forall x \forall y, \quad x + y = (-1)(-x) + (-1)(-y) = (-1)(-x-y) = (-1)(-(y + x)) = y + x$.
\newline (b) If $\alpha \in F$, then $\alpha (0) = 0$ and $ \alpha (-x) = -(\alpha x)$ since $\alpha$ is an endomorphism of $V$. 
\end{rmk}

\begin{defn} (Definition 2.1., \cite{Andre})
Let $(V, F)$ be an $F$-group. The \textit{quasi-kernel} $Q(V)$ (or just $Q$ if there is no danger of confusion) of $(V, F)$ is the set of all $u \in V$ such that, for each pair $\alpha, \beta \in F$, there exists a $\gamma \in F$ for which
\begin{equation}
\alpha u  + \beta u  = \gamma u .
\label{eq:2.1}
\end{equation}
\end{defn}

\begin{lem}{\label{quasiprop}} (\cite{Andre}) 
The quasi-kernel $Q$ has the following properties:
\newline (a) $0 \in Q$;
\newline (b) For $u \in Q \backslash \{0\}$, $\gamma$ in (\ref{eq:2.1}) is uniquely determined by $\alpha$ and $\beta$; 
\newline (c) If $u \in Q$ and $\lambda \in F$, then $\lambda u  \in Q$, i.e. $Fu \subseteq Q$;
\newline (d) If $u \in Q$ and $\lambda_{i} \in F$, $i=1,2,\ldots,n$, then $\sum_{i = 1}^{n} \lambda_{i} u  = \eta u  \in Q$ for some $\eta \in F$;
\newline (e) If $u \in Q \backslash \{0\}$ and $\alpha, \beta \in F$, then there exists a $\gamma \in F$ such that $\alpha u  - \beta u  =\gamma u $.
\end{lem}

\begin{defn}  (Definition 2.3., \cite{Andre})
$(V, F)$ is said to be a \textit{linear} $F$-group if $V = \{0\}$ or $Q(V) \neq \{0\}$.
\end{defn}

\begin{defn}(Definition 4.1., \cite{Andre})
$(V, F)$ is called a \textit{near-vector space over F}  if the following condition holds:
\newline The quasi-kernel $Q = Q(V)$ of $V$ generates the group $(V, +),$ i.e. every element is equal to a finite sum of elements of the quasi-kernel.
\end{defn}

\begin{defn}
We say a near-vector space $(V, F)$ is commutative if for every $\alpha, \beta \in F$ and every $v\in V$, $\alpha(\beta(v))=\beta(\alpha(v))$.
\end{defn}

\begin{defn}
	In a near-vector space $V$ with quasi-kernel $Q$, a \textit{basis} for $Q$ is a minimal generating  subset of $Q$ (note that as $Q$ is not closed under addition, we have that $Q$ is a subset, but not necessarily equal to the span of its basis). By a \textit{basis of $V$} we mean a basis of $Q$ (it is convenient to take the generating elements inside $Q$);  we define $\dim V := \dim Q$. We say a set of elements of $Q$ is independent if it is a non-empty subset of a basis of $Q$.
\end{defn}

\begin{rmk} 
Note that:
\begin{enumerate}
	\item Every near-vector space has a basis.
	\item Every vector space is a near-vector space, with quasi-kernel equal to the whole of $V$, thus the notions of basis coincide here.
\end{enumerate}
\end{rmk}

\begin{defn} (Definition 2.3., \cite{How})
If $(V,F)$ is a near-vector space and $\emptyset \neq V'\subseteq V$ is such that $V'$ is the subgroup of $(V, +)$ generated additively by $FX = \{\lambda x \,|\, x \in X, \lambda \in F\}$, where $X$ is an independent subset of $Q(V)$, then we say that $(V', F)$ is a subspace of $(V, F)$, or simply $V'$ is a subspace of $V$ if $F$ is clear from the context.
\end{defn}

 We can construct near-vector spaces that are not vector spaces from fields and multiplicative automorphisms on them.

 \begin{ex}
Let $F=\mbb F_5$, there is a multiplicative automorphism $\sigma: \mbb F_5 \rightarrow \mbb F_5$ such that $\sigma (a)=a^3$. We can consider the action of $F$ on  $V= \mbb F_5\oplus \mbb F_5$   where the first co-ordinate is acted on using standard multiplication and the second using multiplication by $\sigma(f)$ (we will refer to this as \textit{twisted} multiplication), i.e. for $f\in F$, $f(v_1,v_2)=(fv_1, \sigma(f)v_2)$. So for example $3(2,2)=(1, 4)$. This structure is then a near-vector space, but not a vector space. The quasi-kernel is $(\mbb F_5\oplus 0 )\cup (0\oplus \mbb F_5)$.
\end{ex}

For a near-vector space the action of $F$ on $V$ is, in general, \textit{not} commutative (although for this paper we assume a commutative action). Non-commutative examples can be constructed from near-fields,  see Example 3.1.3 in \cite{karinsthesis}, for a detailed description of this.



\begin{lem}(Lemma 4.5., \cite{Andre}){\label{basis}}
Let $V$ be a near-vector space and let $B = \{u_{i} \,|\, i \in I \}$ be a basis of $Q$. Then each $x \in V$ is a unique linear combination of elements of $B$, i.e. there exists $\xi_{i} \in F$, with $\xi_{i} \neq 0$ for at most a finite number of $i \in I$, which are uniquely determined by $x$ and $B$, such that
\[x = \sum_{i \in I}\xi_{i}u_{i} .\]
\end{lem}

\begin{defn} (Definition 3.2., \cite{HowM2} )
We say that two near-vector spaces $(V_1,F_1)$ and $(V_2,F_2)$ are isomorphic (written $(V_1,F_1) \cong (V_2,F_2))$ if there are group isomorphisms $\theta:(V_1,+) \rightarrow (V_2,+)$ and $\eta:(F_1^{*},\cdot) \rightarrow (F_2^{*},\cdot)$ such that $\theta(\alpha x) = \eta(\alpha) \theta(x)$ for all $x \in V_1$ and $\alpha \in F_1^{*}$.
\end{defn}

\begin{prop}
Let $(V_1,F)$ and $(V_2 , F)$ be near-vector spaces over $F$. Suppose that we have a set-isomorphism (i.e. bijection) between their quasi-kernels $\Phi:Q(V_1)\cong Q(V_2)$, which respects scalar multiplication and addition where it is defined, then $\Phi$ extends to an isomorphism of near-vector spaces.  
\end{prop}

\textit{Proof:} 
First note that if $B\subseteq Q(V_1)$ is a basis for $V_1$ then $\Phi (B)$ will be a basis for $Q(V_2)$, because $\Phi$ is a bijection that preserves scalar multiplication and addition where it is defined. Note also that if $u\in Q(V_1)$ such that $\alpha u+\beta u=\gamma u$	then $\alpha \Phi(u)+\beta \Phi(u)=\gamma \Phi (u)$.

Let $v\in V_1$ then $v=\Sigma_{i=1}^n \lambda _iu_i$ for $\lambda_i\in F$ and $u_i\in B$. Define $\theta(v)=\Sigma_{i=1}^n \lambda _i\Phi(u_i)$. We know from Lemma \ref{quasiprop} that if $u\in Q(V_1)$ then $\lambda u\in Q(V_1)$ for any $\lambda\in F$, so 
$\lambda \Phi(u)=\Phi(\lambda u),$ so $\theta(v)=\Sigma_{i=1}^n \Phi(\lambda _iu_i)$. By the uniqueness of the linear combination (Lemma \ref{basis}) $\theta$ will be onto and one to one. For $v \in V_1$ and $\alpha \in F,$  $\theta(\alpha v) = \alpha \Phi(v)$ follows from the fact that $\alpha$ is an automorphism of $(V_1,+).$ Finally 
$\theta (u_i+u_j)=\Phi(u_i)+\Phi(u_j)=\theta(u_i)+\theta(u_j)$ by definition for all $u_i,u_j\in B$, thus $\theta$ is a near vector space isomorphism (paired with identity on $F$). \\

\QED

\begin{defn} (\cite{Andre}, Section 2)
Let $(V, F)$ be a linear F-group, and let $u \in Q(V) \backslash \{0\}$. Define the operation $+_{u}$ on $\alpha, \beta \in F$ by  
\[(\alpha +_{u} \beta) := \gamma \mbox{ if } \alpha (u)  + \beta (u) =\gamma (u).\]

Note that $\gamma \in F$ will be unique as $u\in Q(V)$. 
\end{defn}

\begin{defn} (Definitions 4.1. and 4.11, \cite{Andre})
We say that $u,v\in Q(V)$ are \textit{compatible} if there is a $\lambda \in F\setminus \{0\}$ such that $u+\lambda v\in Q(V)$. We say that a near-vector space is \textit{regular} if all elements of $Q(V)\setminus \{0\}$ are compatible.
\end{defn}

\begin{lem}(Lemma 4.8., \cite{Andre}){\label{regularitycharacterisation}}
The elements $u$ and $v$ of $Q(V) \backslash \{0\}$ are compatible if and only if there exists a $\lambda \in F \backslash \{0\}$ such that $+_u=+_{\lambda v}$.
\end{lem}

\begin{rmk}
We know from \cite{Andre} (Theorem 4.9.) that compatibility induces an equivalence relation on $Q(V)$, where $u\sim v$ if $+_u=+_{\lambda v}$ for some $\lambda \in F\setminus \{0\}$.
\end{rmk}
\medskip
Moreover, this equivalence relation was used to prove the following result.
\medskip
\begin{thm}(The Decomposition Theorem, Theorem 4.13 from \cite{Andre}\label{decompthm})
\newline Every near-vector space $V$ is the direct sum of regular near-vector spaces $V_{j}$ ($j \in J$) such that each $u \in Q\backslash\{0\}$ lies in precisely one direct summand $V_{j}$. The subspaces $V_{j}$ are maximal regular near-vector spaces.
\end{thm}

\begin{defn}\label{blockdef}
The $V_j$'s obtained from the Decomposition Theorem (Theorem \ref{decompthm}) above are referred to as the \textit{blocks} of $V$. 
\end{defn}

In the case where $F$ is commutative  these blocks have particularly nice properties (see Chapter \ref{chpt3}).
\subsection{Model theory}


Model theory deals with mathematical structures from the point of view of logic. One chooses a language $\mathcal L$ and a set of sentences (or axioms) in that language, called an \textit{$\mathcal L$-theory}, $T$. One can then construct\textit{ $\mathcal L$-structures}. These consist of an underlying set $M$, and an interpretation of the language within that set.  If an axiom $\sigma$ is true under the given interpretation, we write $\mathcal M\models \sigma.$ If for every $\sigma \in T$, $\mathcal M\models \sigma$ then we write $\mathcal M\models T$, and say the $\mathcal L$-structure $\mathcal M$ is a \textit{model} of $T$.  There are several detailed introductions to model theory, we refer the reader to \cite{Hodges}, \cite{Marker}, \cite{kirby} for more details.

Model theorists study the behaviour of definable sets within models of particular theories. A \textit{definable set over $A$}, a set of parameters, is the set of solutions to a formula in $\mathcal L$ using parameters from $A$. 
The language, $\mathcal L$, chosen to describe a particular structure is important, as logical properties of the structure differ according to the language chosen. This is probably easiest to see through an example.

\begin{ex}
There are two ways one could look at vector spaces:
\begin{enumerate}
	\item The standard language chosen for vector spaces over a field $F$ is $\mathcal L=\{+, 0, (f)_{f\in F}\}$, where each $f$ is a unary function interpreted in the structure as scalar multiplication by $f$.
	\item One could also see a vector space as a two sorted structure $(F,V)$, with the field language ($\mathcal L_{\mbox{ring}}=\{+, \times, 0,1\}$) on $F$, the abelian group language ($\mathcal L_{\mbox{ab}}=\{+,0\}$) on $V$, and with a function $\cdot:F\times V\rightarrow V$ representing scalar multiplication.
\end{enumerate} 
From the point of view of model theory these structures are quite different. The latter language is a lot more expressive, as such the definable sets are a lot more complicated. For example, given a vector $v\in V$, $span (v)$ can be defined as the solution set of ``$\exists a \in F (a\cdot v=x)$''. This is not possible in the standard language, where one can find, for any member  $w\in span(v)$  a formula such that $w$ is the unique solution to that formula (i.e. ``$x=\alpha v$'' where $w=\alpha v$), but the whole  of $span(v)$ is not the solution to any formula in the language. One can therefore not express the sentence `every vector is in the span of the set $B$' in this language (it would require infinitely many disjunctions of formulas).
\end{ex}

Generally, the more complicated the definable sets in a structure are, the more difficult it is to analyse it from the point of view of model theory. Hence the `standard' choice of language for vector spaces is less expressive. 

It is extremely useful to be able to reduce the number of quantifiers used to describe definable sets (which may use any number of quantifiers). It is therefore often the first step when analysing a structure (or set of structures) from the model theoretic point of view to establish a language in which the use of quantifiers is redundant.

\begin{defn}
We say an $\mathcal L$-theory $T$ admits \textit{elimination of quantifiers} (or has QE) if for every $\mathcal L$-formula $\phi( x)$ there is a quantifier free $\mathcal L$-formula  $\psi( x)$ such that:
\[T\models \forall  x(\phi( x) \iff \psi( x)).\]

We say $\phi( x)$ and $\psi( x)$ are \textit{equivalent modulo $T$.}

\end{defn}

\begin{ex} (see, for example, \cite{Hodges})
The theory of vector spaces in the `standard' language has QE. Here, all formulas in one variable are equivalent to Boolean combinations of linear equations in one variable. 
\end{ex}


One of the main achievements of geometric stability theory is the classification of theories through combinatorial properties of the definable sets in a theory. Although the classification now goes beyond, in this paper we deal only with stable theories. If one expands $F$ appropriately (see Section \ref{modeltheorynvs}), one can see that a commutative near-vector space is a module (in the standard module language). These are known to be stable (see \cite{Prest}), and in fact we show commutative near-vector spaces are (Theorem \ref{fnvstt}) totally transcendental. The rest of this section gives explanations of the relevant properties.

\begin{defn}
An \textit{$n$-type} of an $\mathcal L$-theory $T$ is a (possibly infinite) set of consistent $\mathcal L$-formulas with at most $n$ free variables. If $\mathcal M\models T$ and $A\subseteq M$, then an \textit{$n$-type over $A$} is one where the formulas are in the language  $\mathcal L_A$ (i.e. $\mathcal L$-formulas with parameters from $A$). Given an $n$-tuple $\bar a \in M^n$ we write $tp(\bar a/A)=\{\phi(\bar x)\in \mathcal L_A: \mathcal M\models \phi(\bar a)\}$.
Sometimes we are lax in our notation and write $\bar a \in M$ or even $a\in M$ to refer to an $n$-tuple of elements of $M$.

An $n$-type (over $A$) is called complete if for every formula $\phi(\bar x)\in \mathcal L$ either it or its negation is in the type.
Suppose $M\models T$ and $A\subseteq M$, then by \textit{$S_n(A)$} we mean the set of complete $n$-types over $A$ in $T$.

\end{defn}

A very useful technique for showing quantifier elimination is the back and forth method. A model, $\mc M$, is \textit{$\omega$-saturated} if all types over finite sets are realised in that model, that is to say that for any finite $A$ and $p\in S_n(A)$ there is a $b\in M$ such that $tp(b/A)=p$.  For quantifier elimination it is then sufficient to show the $\omega$-saturated models have the back and forth property. That is to say:

\begin{thm} \textit{Back and forth, see \cite{Hodges} for a proof.} Let $T$ be an $\mathcal L$-theory, then the following are equivalent:
\begin{enumerate}
	\item $T$ has quantifier elimination.
	\item Let $\mathcal M$ and $\mathcal N$ be $\omega$-saturated models of $T$. Let $\bar a$, $\bar b$ be $n$-tuples from $ M$, $ N$ respectively such that we have a partial isomorphism $\phi: \bar a \mapsto \bar b$. Suppose $c\in M$, then we have a $d\in N$ such that there is a partial isomorphism $\phi': \bar a c \mapsto \bar b d$, vice versa for any $d\in N$ we have a $c\in M$ such that there is a partial isomorphism $\phi': \bar a c \mapsto \bar b d$. In this situation we say the set of partial isomorphisms from $\mathcal M$ to $\mathcal N$ has the back and forth property.
\end{enumerate}
\end{thm}

Quantifier elimination can be used (amongst other things) to understand the types in a particular theory, enabling us to count numbers of types, which is key to understanding where a particular theory is in the stability hierarchy.

\begin{defn}
A theory $T$ is \textit{$\lambda$-stable}, for some infinite cardinal $\lambda$, if for all models $M\models T$ and subsets of $A\subseteq M$ with $|A|<\lambda$, $|S_1(A)|\leq \lambda$.
\end{defn}

We have the following hierarchy within stable theories:

\begin{defn}
\begin{enumerate}
\item A theory is \textit{stable} if it is $\lambda$-stable for some cardinal $\lambda$.
\item A theory, $T$ is \textit{superstable} if there is a cardinal $\kappa$ such that $T$ is $\lambda$-stable for all $\lambda\geq\kappa$. 
\item A theory $T$ is \textit{totally transcendental (t.t.)} if it is $\lambda$-stable for all $\lambda\geq|\mathcal L|$.
\end{enumerate}
\end{defn}

\begin{example}
Vector spaces over a field $F$ are totally transcendental. First note that in this case we have $|\mathcal L|=| F|$ (unless $F$ is finite, in which case $|\mathcal L|=\aleph_0$, see \cite{Marker}).  Given a set of parameters $A$, with $|A|<|\mathcal L|$, we get the following types:
\begin{enumerate}
	\item One type for every element in $span (A)$, so $|span(A)|$-many types.
	\item One type expressing that the vector is not in $span(A)$.
\end{enumerate} 
So the number of types over $A$ is $|span(A)|+1\leq max\{\aleph_0, |F|\}$, thus $|S_1(A)|\leq |\mathcal L|$, so the theory is totally transcendental.
\end{example}

Totally transcendental theories are very well understood, and carry a good notion of dimension. We introduce Morley rank and state some properties of Morley rank in totally transcendental theories. The reader is referred to \cite{pillaystablelecturenotes} for more details (see also Section 2.5 in \cite{Hodges} for the definition of an elementary extension).

\begin{defn}
Let $T$ be a complete theory, $\mc M\models T$. We first define, by induction, what it means for a formula, $\phi(\bar x, \bar a)$ (where $\bar a$ are parameters from $M$), to have Morley rank greater than or equal to some ordinal $\alpha$ (written $RM(\phi(\bar x, \bar a))\geq \alpha$).  
\renewcommand{\labelenumii}{(\alph{enumii})}
\begin{enumerate}
\item  $RM(\phi(\bar x, \bar a))=-1$ iff $\mc M\models \neg\exists \bar x \phi(\bar x, \bar a)$.
\item  $RM(\phi(\bar x, \bar a))\geq 0$ if $\mc M\models \exists \bar x \phi(\bar x, \bar a)$.
\item $RM(\phi(\bar x, \bar a))\geq \alpha +1$ if there is some elementary extension $\mc N$ of $\mc M$
and there are formulas $\psi_j(\bar x,\bar b_j)$, $j\in\omega$ with $\bar b_j\in N$, such that:
\begin{enumerate}
\item $\mc N\models \psi_j(\bar x,\bar b_j)\rightarrow \phi(\bar x, \bar a)$ for all $j\in \omega$.
\item $RM(\psi_j(\bar x,\bar b_j))\geq \alpha$ for all $j\in\omega$.
\item The sets that the $\psi_j(\bar x,\bar b_j)$ define are pairwise disjoint (i.e.for $i\not =j$ $\mc N\models \neg \exists \bar x(\psi_j(\bar x,\bar b_j)\wedge \psi_i(\bar x,\bar b_i))$).
\end{enumerate}

\item  For $\alpha$ a limit ordinal $RM(\phi(\bar x, \bar a))\geq \alpha$ if $RM(\phi(\bar x, \bar a))\geq \delta$ for all $\delta<\alpha.$

\end{enumerate}

A formula $\phi(\bar x,\bar a)$ has \textit{Morley rank} $\alpha$ (written $RM(\phi(\bar x,\bar a))=\alpha$) if $\alpha$ is the largest ordinal such that $RM(\phi(\bar x,\bar a))\geq \alpha$. If there is no such $\alpha$ we say $RM(\phi(\bar x,\bar a))=\infty$.

\end{defn}

\begin{prop}\label{dMprop}
Let $RM(\phi(\bar x, \bar a))=\alpha$, then there is a greatest integer $d$ such
that there exist $\psi_j(\bar x, \bar b_j)$ for $1\leq j \leq d$ such that $RM( \psi_j(\bar x,\bar b_j)) = \alpha$ for $1\leq j\leq k$ and $M \models\psi_j(\bar x, \bar b_j)\rightarrow \phi (\bar x, \bar a)$ for $1\leq j \leq d$, and the $\psi_j(\bar x, \bar b_j)$ are pairwise disjoint.
\end{prop}

\begin{defn}
The \textit{Morley degree} of a formula $\phi(\bar x, \bar a)$ with ordinal valued Morley rank is precisely the $d$ obtained from Proposition \ref{dMprop}. We write $dM(\phi(\bar x,\bar a))=d$.
\end{defn}

\textit{Notation:} For a definable set $X$, we will often use $RM(X)$ and $dM(X)$ to refer to the Morley rank and degree of the formula defining $X$. For a type $tp(a/B)$ we write $RM(a/B)$ for the minimum Morley rank of the formulas contained in the type. The \textit{Morley rank of a structure} $M\models T$ is $RM(x=x)$. If this is ordinal valued, then its Morley degree is $dM(x=x)$.

\begin{rmk}
 A theory $T$ is totally transcendental if all formulas in the theory have ordinal valued Morley rank. That is to say for every definable set $X$, $RM(X)<\infty$.
\end{rmk}

\begin{ex}
Vector spaces have minimal Morley rank and degree, i.e. in a vector space $RM(x=x)=dM(x=x)=1$. In this case we can use quantifier elimination to see that every definable set in one variable is either finite or cofinite, the result is immediate from this. 
\end{ex}

In fact structures with Morley rank and degree equal to $1$ are called \textit{strongly minimal} and are the most basic in the geometric model theory hierarchy. As well as vector spaces, strongly minimal structures include algebraically closed fields (in the field language) and infinite sets in the empty language (note that by convention all languages contain identity). 
 Morley rank in algebraically closed fields is equivalent to transcendence degree. In vector spaces the independence relation given by Morley rank (i.e. $a$ is independent from $B$ over $C$ if $RM(a/B\cup C)=RM(a/C)$) is linear independence. See \cite{Marker} for more details on strongly minimal structures.


\section{Commutative near-vector spaces}\label{chpt3}

In this section we consider a near-vector space $(V,F)$ where we will assume $(F,\circ)$ is commutative (i.e. for every $\alpha,\beta \in F$ and $v\in V$ we have $\alpha(\beta(v))=\beta(\alpha(v))$). We establish some basic facts, and give a more precise decomposition theorem in this setting. Recall that $B$ was a basis for the near-vector space, and that $B\subseteq Q(V)$. 

\begin{rmk}
$(F,+_u,\circ)$ is a field.
\end{rmk}

\textit{Proof:}
By (\cite{karinsthesis} Theorem 2.3.5) this must be a near field, that it is a field follows from commutativity of $+_u$ (obvious), and commutativity of $\circ$ (assumed). \\

\QED

\begin{Fact}
From Lemma \ref{basis} and the remark above, it is not difficult to see that $V\cong \oplus_{u\in B} (F, +_u)$  (with $B$ a basis for $V$). 
\end{Fact}

\begin{lem} {\label{regularitycharacterisation2}}If $(F, \circ)$ is commutative then $+_u=+_{\lambda u}$ for all $\lambda\in F$.
\end{lem}
\textit{Proof}
Let $\alpha,\beta \in F$, then:

\[\begin{array}{lllllll}
\lambda((\alpha +_{\lambda u} \beta)(u))&=& 	(\alpha +_{\lambda u} \beta)(\lambda(u))& F\mbox{ is commutative}\\
&=& 	\alpha(\lambda(u)) + \beta(\lambda(u))&\\
&=& 	\lambda(\alpha(u)) + \lambda(\beta(u))& F\mbox{ is commutative}\\
&=& 	(\lambda\alpha +_u \lambda\beta)(u)\\
&=&   (\lambda(\alpha +_u \beta))(u)& (F, +_u, \circ) \mbox{ is a field} \\
&=& 	\lambda(\alpha +_u \beta)(u)& \\
\end{array}\]

We can precompose by $\lambda^{-1}$ to get $(\alpha +_{\lambda u} \beta)(u)=(\alpha +_u \beta)(u)$, so as $F$ is fixed point free $\alpha +_{\lambda u} \beta = \alpha +_u \beta$, i.e.  $ +_{\lambda u} =+_u $.
\QED

Using the compatibility equivalence relation we can partition $Q \backslash \{0\}$ into sets $Q_{j}$ ($j \in J$) of mutually pairwise compatible vectors. Furthermore, let $B \subseteq Q\backslash\{0\}$ be a basis of $V$ and let $B_{j} := B \cap Q_{j}$. By our partitioning, the $B_{j}$'s are disjoint and each $B_{j}$ is an independent subset of $B$.
Let $V_{j}$ := $\langle B_{j} \rangle$ be the subspace of $V$ generated by $B_{j}$, these are the {\it{blocks}} of $V$ as defined in Definition \ref{blockdef}. By the Decomposition Theorem each block is a maximal regular near-vector space over $(F, \circ)$. But since $(F, \circ)$ is commutative, in this case we have more:

\begin{thm}\label{com+reg=vs}
A near-vector space $(V,F)$ over commutative $F$ is regular if and only if it is a vector space over $(F,+_u,\circ)$ for any $u\in Q(V)\backslash\{0\}$.
\end{thm}

\textit{Proof:} Clearly every vector space over a field is regular as $Q(V)=V$. Conversely, 
suppose $u,v\in Q(V)\backslash\{0\}$, then as $V$ is regular, by Lemma \ref{regularitycharacterisation} we have that there is a $\lambda \in F\setminus \{0\}$ such that $+_{\lambda u}=+_v$ and by  Lemma \ref{regularitycharacterisation2} we have that $ +_{\lambda u} =+_u,$ so that $+_u = +_v.$ Now put $+ = +_u,$ then $Q(V) = V$ and since $(F,+,\circ)$ is a field, $V$ is a vector space over $(F,+,\circ).$ 
\QED

\begin{thm}
Each block $V_{j}$ := $\langle B_{j} \rangle$ is a vector space over $(F, +_u,\circ)$ for $u \in Q(V_j)\backslash\{0\}.$ 
\end{thm}
\textit{Proof:}
Since each $V_{j}$ is regular, the result follows from the above theorem.
\QED

The next result immediately follows from the previous result and Theorem \ref{decompthm}.

\begin{cor} \label{blocksareVS}
 Every commutative near-vector space $V$ is isomorphic to the direct sum of blocks where each block is a vector space over $(F, +_b, \circ)$ for $b\in B$. The dimension of the block (as a vector space over $(F, +_b, \circ)$) is the number of elements of $B$ equivalent to $b$. 
\end{cor}

\begin{ex}
If $F$ is finite and commutative all of the $(F, +_u, \circ),$ for  $u \in Q(V)\backslash\{0\},$ are finite fields of the same size, so isomorphic. 
So the only near-vector spaces over finite fields are those where each basis element is ``twisted'' by a field automorphism. Note that if they are all twisted by the same automorphism then this becomes a vector space over the finite field.  This gives us a quick way of finding all near-vector spaces over a finite field through the automorphism group of the field. See \cite{HowM2} for more results on near-vector spaces over finite fields.
\end{ex}

\begin{ex}\label{QF3}

If $F$ is infinite then we no longer have that all the $(F, +_u, \circ),$  for $u \in Q(V)\backslash\{0\},$ need to be isomorphic as fields. However the underlying multiplicative groups must be isomorphic. 

For example, if $\mbb F_3$ is the field with three elements, $\mbb F_3(t)$ is the quotient field and $\mathbb F_3[t]$  the polynomial field of $\mbb F_3$ in the indeterminate $t$, then $\mathbb Q\oplus \mathbb F_3(t)$ is a near-vector space over $\mathbb Q.$ 

First consider $ \mathbb Z$ which is generated by primes. Now, $\mathbb F_3[t]$ is a PID, so generated by the generators of countably many principal ideals. Therefore primes in $\mbb Z$ can be sent to generators of these principal ideals (and $1,-1$ to $1,-1$ respectively). Then take the induced map on the fraction fields. Thus $(\mathbb Q\setminus \{0\},\cdot)\cong (\mathbb F_3(t)\setminus \{0\},\cdot)$,  call this isomorphism $\sigma$.

It is then clear that $\mathbb Q\oplus \mathbb F_3(t)$ is a near-vector space because $\sigma:(\mathbb Q\setminus \{0\},\cdot)\cong (\mathbb F_3(t)\setminus \{0\},\cdot)$, so $\mbb Q$ can be seen as a subset of $Aut(\mathbb F_3(t))$, acting via this isomorphism (i.e. for $\lambda\in \mbb Q$, $v\in \mathbb F_3(t)$ define $\lambda(v)=\sigma(\lambda)\cdot v$) .


\end{ex}


The following theorem was given in \cite{Andre} giving necessary and sufficient conditions for a near-vector space to be a vector space.

\begin{thm} (Theorem 5.5 from \cite{Andre}) \label{andrescharacterisationofvs}
Let $(V,F)$ be a near-vector space with $\mbox{dim}(V) > 1$. $F$ is a division ring and $V$ a vector space over it if and only if $V$ coincides with its quasi-kernel $Q(V)$.
\end{thm}

However, \cite{How} (Example 4.2) notes the following counter-example:

\begin{ex}
Let $F=\mathbb R$, and $V=\mathbb R\oplus \mathbb R$ where for $\alpha\in F$, and $(v_1,v_2)\in V$, $\alpha(v_1,v_2)=(\alpha ^3 v_1, \alpha^3v_2)$. It can be verified that this is a near-vector space, and we have that $Q(V)=V$. However, $V$ is not a vector space over $(\mbb R, \cdot, +)$ with the usual addition, as it is not distributive. Suppose $\alpha,\beta \in F\setminus \{0\}$, then 
\[\begin{array}{lllllll}
	(\alpha+\beta)(1,1)&=&((\alpha+\beta)^3, (\alpha+\beta)^3)\\
&\ne & (\alpha^3+\beta^3,\alpha^3+\beta^3)\\
&=& \alpha(1,1)+\beta(1,1). 
\end{array}\]

However, this counterexample is still in fact a vector space, it is a vector space over the field $(F, \circ, +_u)$, where for $\alpha,\beta \in F$, we have $\alpha +_u\beta= (\alpha^3+\beta^3)^{\frac{1}{3}}$. This field is in fact isomorphic to $(\mbb R, \cdot, +)$ via the map taking $\alpha$ to $\alpha^3$.

Interestingly this near-vector space is isomorphic as a near-vector space to the two dimensional vector space over $\mathbb R$, but as vector spaces they are not isomorphic as vector space isomorphisms assume you are working over the same field. 
\end{ex}

We conclude that Theorem \ref{andrescharacterisationofvs} is correct, but that the field over which the near-vector space is a vector space needs some clarification. As in the definition of near-vector spaces $F$ is not equipped with an addition, it is not immediately obvious what field is being referred to in Theorem  \ref{andrescharacterisationofvs}. The field should be $(F, \circ, +_u)$, once this is established the statement holds.




\section{Model theory of near-vector spaces}\label{modeltheorynvs}

We wish to study $F$-groups and near-vector spaces from a model theoretic perspective. We start with some basic notation and results. We then show that there is no first-order theory whose models are precisely the near-vector spaces. We then go on to show that all models of the theory of a near-vector space with finitely many blocks are near-vector spaces with the same block type. We show that this complete theory has quantifier elimination and is totally transcendental, of Morley rank the number of blocks.

For ease of understanding we have assumed $F$ to be commutative throughout, although we suspect some of the work would carry through to the non-commutative case.

\subsection{Notation}

Here we fix a commutative $F\subseteq Aut(V)$, we use the language \\$\mc L_{Fnvs}=\{+, 0, (\lambda)_{\lambda\in F}\}$ where each $\lambda$ is seen as a unary function symbol.
Given any near-vector space $V$ over $F$ we can see $V$ as an $\mc L_{Fnvs}$-structure by interpreting $+$ as addition and $\lambda$ as scalar multiplication. 

\begin{rmk}
We can express the axioms of an $F$-group using infinite axiom schemes in this language. 
We can quite easily  also add (infinitely many) axioms that express commutativity of $F$ or that $V$ is a linear $F$-group.\end{rmk}

In this language, however, we cannot add first-order sentences which express that an $F$-group is generated by elements of the quasi-kernel. This would involve an infinite number of disjunctions. In fact the notion of a near-vector space is not first-order in any language, we can see this through the following example.

\begin{ex}
Let $V=\oplus_{i\in \mbb N} \mbb R$, with the standard co-ordinate-wise addition, and let $F=\mbb R$ and $\alpha\in F$ act on $v=( v_i) _{i\in \omega}$ as follows:

\[ \alpha (v)=(\alpha^{2i+1} v_i)_{i\in \omega}\]

This is clearly an $F$-group. The quasi-kernel is the set of vectors with at most one non-zero vector, this clearly generates the whole space (recall almost all the co-ordinates of $V$ are $0$). So $V$ is a near-vector space. However, $V$ has an ultrapower that is not a near-vector space, showing that this concept is not first-order.

Let $B=\{b_1,...,b_n,....\}\subset Q(V)$ be a set that generates the near-vector space, we can assume that $b_i=(v_j)_{j\in \mbb N}$ with $v_i=0$ for $i\ne j$.  Each $b_i$ generates a block $B_i$ of $V=\oplus_{i\in \mbb N} B_i$.

Let $\mathcal W=\prod V/\mathcal U$, with $\mathcal U$ any non-principal ultrafilter on $\omega$.
 Consider the following element of this ultraproduct:

\[w=(w_j)/\mathcal U \mbox{ where } w_j= \sum_{i=1}^j b_i\]

Now each $w_j$ is in the span of exactly $j$ elements of $Q(V)$, but is not in the span of any fewer elements of $Q(V)$. 

Notice that the quasi-kernel $Q(\mathcal W)$ is defined to be the set of elements of $ u\in  \mathcal W$ such that for every $\alpha, \beta \in \mbb R$ there is a $\gamma \in \mbb R$ such that $\alpha (u)+\beta(u)=\gamma(u)$. 
Suppose $u=(u_j)/\mathcal U\in Q(\mathcal W)$ where each $u_j\in V$, then $\alpha u=(\alpha u_j)/\mathcal U$. Suppose for $\alpha, \beta\in F$,  $\alpha (u)+\beta(u)=\gamma(u)$, then for all but finitely many $u_j$ we have $\alpha (u_j)+\beta(u_j)=\gamma(u_j)$. As this holds for any $\alpha, \beta \in F$, almost all the $u_j$'s must be in the same block (or zero).

So for any $n\in \mbb N$ we cannot express $w$ as a sum of $n$ elements of the quasi-kernel, as we would need more than $n$ elements of $Q(\mathcal W)$ to express $w_{j}$ for $j>n$. 
Therefore in the ultrapower $w$ will not be in the span of the quasi-kernel of $\mathcal W$, so this is not a near-vector space.

Now consider the element $t=(u_j)/\mathcal U\in \mathcal W$ such that $u_j=(v_j^i)_{i\in \omega}\in V$ with:

\[v^i_j=\left\{\begin{array}{llll}
	1&\mbox{ if } i=j\\
	0&\mbox{ if } i\ne j\\	
\end{array}\right.\]
Now for each $j$ we have that $u_j\in Q(V)$, thus $t\in \prod Q(V)/\mathcal U$. However, for a given $\alpha, \beta \in F$ we have that $\alpha (u_j)+\beta (u_j)= (\alpha^{2j+1}+\beta^{2j+1})^{\frac{1}{2j+1}}(u_j)=\gamma_j(u_j)$, now $t\notin Q(\mathcal W)$ as 
$$(\alpha^{2j+1}+\beta^{2j+1})^{\frac{1}{2j+1}}\ne(\alpha^{2k+1}+\beta^{2k+1})^{\frac{1}{2k+1}} \mbox{ for } k\ne j.$$

We therefore have that $Q(\mathcal W)\subsetneq \prod Q(V)/\mathcal U$. The quasi-kernel is therefore not definable in this case.
\end{ex}

We will see later that it is possible to axiomatise certain near-vector spaces, which we call finite block near-vector spaces. To do this we need to introduce some new concepts.

\begin{defn}
In this setting we let $\bar F$ be the expansion of $F$ by formal finite sums of elements of $F$. For $\alpha_1,...,\alpha_n \in F$, we will denote the sums as $\alpha_1 +_.\alpha_2$, or ${\Sigma^n_i}^. \alpha_i$ (the dots to distinguish them from other sums used). Note that $\bar F$ is closed under composition as $F$ is. 
\end{defn}

\begin{defn}\label{Fbaracts}
Suppose $V$ is a near-vector space over $F$. We can see $\bar F$ as acting on $V$ by pointwise addition. That is to say for 
${\Sigma^n_i}^. \alpha_i\in \bar F$, $v\in V$ define 
$( {\Sigma^n_i}^. \alpha_i)(v)=\alpha_1(v)+...+\alpha_n(v)$ 
(where $+$ is the addition in $V$). It is clear that as $\bar F\subseteq End(V)$ actions formed by pointwise addition will also be in $End(V)$. 
\end{defn}

Given a near-vector space $(V,F)$ some elements of $\bar F$ will act as automorphisms on $V$. We choose not to add formal inverses to $\bar F$, as whether elements of $\bar F$ are elements of $Aut(V)$ depends on which $V$ is being acted on. Instead we define $\mbox{frac}_V(\bar F)$ to be the subset of $End(V)$ generated by $\bar F$ and $\bar F\cap Aut (V)$ (that is to say closing under inverses and their sums when they exist). 

It is clear that we can now see $\bar F$ as a subring of $End(V)$, and that if $(V,F)$ is a near-vector space, $F$ commutative, then $V$ is an  $\bar F$-module. For two near-vector spaces over $F$, $\bar F$ will be the same in both cases (as it is just formal sums), but $\bar F$ may well act quite differently depending on the near-vector space. 

\begin{example}
Let $V_1=\mbb Q\oplus\mbb F_3(t)$ and $V_2=\mbb Q.$ By Example \ref{QF3} we can see both of these as near-vector spaces over $F=\mbb Q$. In the case of $V_2$, as $\mbb Q$ is a $\mbb Q$-vector space, using the action defined in Definition \ref{Fbaracts}, $\bar F$ acts exactly as $F$ does. That is to say for every formal sum $\bar \alpha \in\bar F\setminus \{0\}$ there is a $q\in \mbb Q$ such that $\forall v\in V_2$, $\bar \alpha (v)=q(v)$, so $\bar F$ with this action is a subset of $Aut(V_2) \cup \{0\}$.

This is not the case for $V_1$. For example, if we let $\bar \alpha=1+_.1+_.1$ and $(u_1,u_2)\in V_1$ we get 
\[\bar \alpha (u_1,u_2)= ((1+1+1)u_1, (1+1+1)u_2)=(3u_1,0).\]
 So $\bar \alpha$ acts as $3$ on the first co-ordinate and $0$ on the second, it is neither an automorphism of $V_1$, nor $0$, so it cannot be contained in $F$. Therefore $\bar F\not \subseteq Aut(V_2)\cup \{0\}$.

\end{example}

\begin{defn}\label{blocktype}
Let $V$ be a near-vector space over $F$. Then by Corollary \ref{blocksareVS} there is a set $I$ and $B_i$  for $i\in I$ such that each $B_i$ is a vector space over the field $(F,\circ, +_{u_i})$ and  $V=\oplus B_i$. We define \textit{the block type of $V_1$}, as the set $\{+_{u_i}\}_{i\in I}$ and denote this by $BT(V)$. Note that by Theorem \ref{com+reg=vs} this does not depend on the choice of $u_i$.
\end{defn}

\begin{rem}
Suppose $V_1$ and $V_2$ are both near-vector spaces over $F$, and suppose $V_1=\oplus_{i\in I} B_i$ and $V_2=\oplus_{j\in J} C_j$ where $B_i$ is a vector space over $(F,\circ, +_{u_i})$ and $C_j$ is a vector space over $(F,\circ, +_{v_j})$ (i.e. $B_i\cong (F,+_{u_i})^{n_i}$ and $C_j\cong (F,+_{v_j})^{m_j}$). Then $V_1$ and $V_2$ have the same block type if they contain the same type of blocks, i.e. $\{(F,\circ, +_{u_i})\}_{i\in I}=\{(F,\circ, +_{v_j})\}_{j\in J}$ [note here that it is important to have equality, not merely isomorphism].
\end{rem}

\begin{defn}\label{blocktypefinite}
We say that a near-vector space $V$ has finite block type if $BT(V)$ is finite. Equivalently they contain finitely many blocks in the sense of Theorem \ref{decompthm}.
\end{defn}

First we establish some results around conditions under which near-vector spaces are in fact vector spaces. 
Recall that we are assuming $F$ to be commutative. We can expand the language to $L_{\bar{F}nvs}=\{+, 0, (\lambda)_{\lambda\in \bar{ F}}\}$ with the obvious interpretation. Note that as each unary function $\alpha\in \bar F$ is quantifier free definable in $L_{\bar{F}nvs}$, this expansion will not affect any quantifier elimination result. 

Note that if we were to expand the language to include unary predicates for all elements of $frac_V(\bar F)$ then this would a priori have an effect on quantifier elimination, as to define inverse elements of $\bar F\cup Aut(V)$ one would need quantifiers. However, we will see (Lemma \ref{F=frafF}) that in the finite block case $frac_V(\bar F)=\bar F$.

\begin{prop}
Let $F$ be commutative. Then $(V,F)$ is a vector space (i.e. $V$ is a vector space over $F$) if and only if $\bar{F}$ (as a subset of $End(V),$ see Definition \ref{Fbaracts}) is equal to $F$.
\end{prop}

\textit{Proof:}   If $(V,F)$ is a vector space then $F$ is closed under finite sums, so $\bar{F}=F$. Conversely, by Theorem \ref{com+reg=vs} it is sufficient to show that for all $u,v\in Q(V)$, $+_u=+_v$. 
 Suppose $F=\bar{F}$, let  $u,v\in Q(V)$, $\alpha,\beta \in F$ then:
 \[\begin{array}{lllll}
    (\alpha+_u\beta)(u)&=&\alpha(u)+\beta(u)\\
    &=& (\alpha+_\cdot\beta)(u)\\
    &=& (\gamma)(u)& \mbox{ for some }\gamma\in F \mbox{ as }F=\bar{F}.\\
   \end{array}
\]
So as $F$ acts fixed point free, we have $\alpha+_u\beta=\alpha+_\cdot \beta=\gamma$, this is true for all $u\in Q(V)$, so in particular $+_u=+_\cdot=+_v$.\QED


\begin{prop}
 Let $(V, F)$ be a commutative near-vector space, then \newline  $(V,\mbox{frac}_V(\bar{F}))$ is a vector space if and only if $\bar{F}\subseteq Aut(V)\cup \{0\}$.
\end{prop}

\textit{Proof:} 
Clearly if $(V,\mbox{frac}_V(\bar{F}))$ is a vector space, then $\bar{F}\subseteq Aut(V)\cup \{0\}$. Conversely, as we know $V$ is an $\bar{F}$-module it is sufficient to show that $(\mbox{frac}_V(\bar{F}), \circ, +_\cdot)$ is a field. 
\begin{itemize}
 \item Now $(\mbox{frac}_V(\bar{F}), +_\cdot)$ is a commutative group, with $0\in F$ as its identity element. 
 \item $\mbox{frac}_V(\bar{F})$ is closed under composition, has an identity element $1$, and commutativity follows from commutativity of $(F,\circ)$ and $(V,+)$. So we need to show every element of $\alpha\in \mbox{frac}_V(\bar{F})$ has an inverse inside $\mbox{frac}_V(\bar{F})$. As we have assumed $\bar{F}\subseteq Aut(V)$ the result follows from the fact that $\mbox{frac}_V(\bar{F})$ is closed under inverses. 

\end{itemize}



\color{black}{

From a model theory point of view we would not distinguish between $V$ being a vector space over $F$ or $V$ being a vector space over $\mbox{frac}_V(\bar{F})$, as the latter is just a definitional expansion of the former.

\begin{ex}\label{fgnotnvs}
Certainly $\mbox{frac}_V(\bar{F})$ is not always a field, if we consider Example \ref{QF3} then for $\alpha \in \mbb Q\setminus \{0\}$ we have:
\[\alpha +_\cdot \alpha +_\cdot\alpha (v_1,v_2)=(3\alpha v_1, 3\sigma(\alpha)v_2)=(3\alpha v_1,0).\]
This is clearly not an automorphism (take $v_2\ne 0$). \\





\end{ex}

Even in the more general case of finite block near-vector spaces, we still get that $V$ can be recovered from $\bar{F} \cap Aut(V)$. That is to say, if for two near-vector spaces (over $F$) the same elements 
 of $\bar{F}$ act as automorphisms, then the block types of both vector spaces will be the same. This is what we mean when we say that  the theory is `controlled' by which elements of the pointwise additive closure of $F$ are automorphisms of the near-vector space. This is particularly relevant from the point of view of model theory as we can express whether an element of $\bar{F}$ is an automorphism using first-order logic.

The following is a useful lemma, and holds for all commutative near-vector spaces, including those with infinite block type. Essentially it says that additions in the block type of a near-vector space are determined by which sums of elements of $F$ are zero.
\begin{lem}\label{autodetby0}
Let $(V, F)$ be a commutative near-vector space. 
Let $+_u$ and $+_v$ be additions generated by elements $u, v\in Q(V)$, so both $(F,\circ, +_u)$ and $(F,\circ, +_v)$ are fields.  Suppose that for any $\alpha_1, \alpha_2,\alpha_3\in F$ we have:
\[\alpha_1+_{u}\alpha_2+_{u}\alpha_3=0 \mbox{ if and only  if }\alpha_1+_{v}\alpha_2+_{v}\alpha_3=0.\]
Then we have that $+_{u}$ and $+_{v}$ are equal. Moreover, as the converse is clearly true, which sums are zero determine the addition completely.
\end{lem}

\textit{Proof:} Suppose for contradiction that $\alpha+_{u}\beta\ne \alpha+_{v}\beta$ for some $\alpha,\beta \in F$.
Now let $\alpha+_{u}\beta=\gamma_u\in F$, $\alpha+_{v}\beta=\gamma_v\in F$, so $\gamma_u\ne\gamma_v$. 
Clearly, $\alpha+_{u}\beta+_{u}\gamma_u(-1)=0$, but $\alpha+_{v}\beta+_{v}\gamma_u(-1)\ne 0,$ which contradicts our assumption.\QED

We now establish that the block type of a near-vector space $V$ over $F$ completely determines which elements of $\bar{F}$ are automorphisms of $V$. Here we do not have any conditions on the size of $BT(V)$.

\begin{prop}\label{BTfixesauto}

Suppose $V_1$ and $V_2$ are near-vector spaces over $F$, and suppose further that $BT(V_1)=BT(V_2)$, then $\bar{F} \cap Aut{V_1}=\bar{F} \cap Aut{V_2}$.

\end{prop}

\textit{Proof:} 
Throughout this proof we use the notation of Definition \ref{blocktype}. That is to say that $V_1=\oplus B_i$ and $V_2=\oplus C_j$ and $BT(V_1)=\{+_{u_i}:i\in I\}$, $BT(V_1)=\{+_{v_j}:j\in J\}$.
For any $v \in V_1$ we have uniquely determined $v_i'$ for $i\in I$ such that $v=\sum_{i\in I}v_i'$ and $v_i'$ is zero in all summands apart from $B_i$. 
Let $\alpha_1 +_{\cdot}...+_{\cdot}\alpha_n\in F$, then, for $v\in V_1$ we have that $(\alpha_1 +_{\cdot}...+_{\cdot}\alpha_n)(v)= \Sigma _i (\alpha_1+_{u_i}...+_{u_i}\alpha_n )v_i$. 
Now $\alpha_1 +_{\cdot}...+_{\cdot}\alpha_n$ is an automorphism in $V_1$ if and only if every $(\alpha_1+_{u_i}...+_{u_i}\alpha_n )$ is an automorphism (i.e. co-ordinate-wise $\alpha_1 +_{\cdot}...+_{\cdot}\alpha_n$ acts as an automorphism). 
As $(\alpha_1+_{u_i}...+_{u_i}\alpha_n)\in F$, and $F\setminus \{0\}\subseteq Aut(V_1)$, we have that $\alpha_1 +_{\cdot}...+_{\cdot}\alpha_n$ is an automorphism if and only if every $(\alpha_1+_{u_i}...+_{u_i}\alpha_n )\ne 0$. 
Similarly, in $V_2$,  $\alpha_1 +_{\cdot}...+_{\cdot}\alpha_n\in Aut(V_2)$ if and only if for every $j$, $(\alpha_1+_{v_j}...+_{v_j}\alpha_n )\ne 0$.\\

Suppose for contradiction that $\bar{F} \cap Aut{V_1}\ne\bar{F} \cap Aut{V_2}$, then without loss of generality there are some $\alpha_1,...,\alpha_n\in F$ such that $\alpha_1 +_{\cdot}...+_{\cdot}\alpha_n\in Aut (V_1)$ and 
$\alpha_1 +_{\cdot}...+_{\cdot}\alpha_n\notin Aut (V_2)$, so for all $i$, $(\alpha_1+_{u_i}...+_{u_i}\alpha_n )\ne 0$, but for some $j$, $(\alpha_1+_{v_j}...+_{v_j}\alpha_n )= 0$. Therefore, by the claim $+_{v_j}\ne +_{u_i}$ for all $i$, thus $V_1$ and $V_2$ must have different block type.


\QED

The converse of Proposition \ref{BTfixesauto} is not true in general for all near-vector spaces. Below we give an example of two infinite block type near-vector spaces in which the same elements of $\bar F$ are automorphisms, but whose block types are different.

\begin{ex}
Let $P$ denote the set of all primes, and $S_P$  the set of all permutations on $P$. If $\sigma\in S_P$ then there is a unique extension, $\hat\sigma: \mbb Q\ra \mbb Q$ of $\sigma$ to an automorphism of the multiplicative structure of $\mbb Q$. 

Let  $\mbb Q_\sigma= (\mbb Q, \cdot, +_\sigma)$ be the field twisted by $\sigma$. That is to say for $a,b\in \mbb Q$, $a+_\sigma b=\sigma^{-1}(\sigma(a)+\sigma(b))$ where $+$ denotes the standard addition in $\mbb Q$. 
Now consider $V_1=\bigoplus_{\sigma \in S_P} \mbb Q_\sigma$ and $V_2=\bigoplus_{\sigma \in S_P\setminus\{id\}} \mbb Q_\sigma$. These are both near-vector spaces over $\mbb Q$. Note that if we take $u$ in a summand $\mbb Q_\sigma$ then $u$ will be in the quasi-kernel and $+_u=+_\sigma$. 
It is also clear that $+_{id} \not \in \{+_{\sigma}\}_{\sigma \in S_P\setminus\{id\}}$ , so $BT(V_1)\ne BT(V_2)$.\\\\
\textit{Claim:} $\bar F\cap Aut( V_1)=\bar F \cap Aut(V_2)$.

\textit{Proof of claim:} It is sufficient to show $\bar F\setminus Aut(V_1)=\bar F\setminus Aut(V_2)$. Clearly  $Aut(V_1)\cap \bar F\subseteq Aut(V_2)\cap \bar F$, so $\bar F\setminus Aut(V_2)\subseteq\bar F\setminus Aut(V_1)$ . Let $\alpha={\Sigma_{i=1}^n}{\cdot} \alpha_i \in \bar F$. Suppose $\alpha \notin Aut(V_1)$, then, by reasoning similar to Proposition \ref{BTfixesauto} there is some $\sigma \in S_P$ such that $\alpha_1+_\sigma...+_\sigma \alpha_n=0$. If $\sigma\in S_P\setminus\{id\}$ then $\alpha \notin Aut(V_2)$. Otherwise, if $\sigma =id$ then there is a $\tau \in S_P\setminus \{id\} $ such that  $\alpha_1+_\tau...+_\tau \alpha_n=\alpha_1+_\sigma...+_\sigma \alpha_n=0$. This is because $\alpha_1,...,\alpha_n$ are multiplicatively generated by finitely many primes, so we can choose a $\tau \ne id$ that acts as identity on these primes, so 
\[\begin{array}{lllll}
   \alpha_1+_\tau...+_\tau \alpha_n &=&\alpha_1+...+\alpha_n\\
   &=&\alpha_1+_\sigma...+_\sigma \alpha_n.
  \end{array}
\]
Therefore $\alpha \notin Aut(V_2)$.  \QED
\end{ex}

Note that similar examples of this phenomenon can be constructed using infinite block type near-vector spaces whose block types differ only by \textit{finitely} many elements.

The converse of Proposition \ref{BTfixesauto} \textit{is} true in the context of finite block near-vector spaces. We establish this below.

\begin{prop}\label{fixedbyauto}
Let $V_1$ and $V_2$ be finite block near-vector spaces over $F$, then $BT(V_1)=BT(V_2)$ if and only if $\bar{F} \cap Aut(V_1)=\bar{F} \cap Aut(V_2)$.
\end{prop}

\textit{Proof:} 
The left to right direction is Proposition \ref{BTfixesauto}.  Let $BT(V_1)=\{+_i:i\in I\}$, $BT(V_2)=\{+_j:j\in J\}$, with both $|I|$ and $|J|$ finite. 
Note that for any near-vector space $V$ over $F$ we have that $\bar F$ is a subring of $End(V)$. Now for each $+_i\in BT(V_1)$ we can define a ring homomorphism:
\[\begin{array}{cccclllll}
	\Phi_i: (\bar F , \circ, +_\cdot) &\rightarrow& (F, \circ ,+_i)\\
	\Phi_i (\alpha_1+_\cdot...+_\cdot \alpha_n)&=& \alpha_1+_i...+_i \alpha_n.
\end{array}\]

As $Im(\Phi_i)=(F, \circ ,+_i)$ is a field, we must have that $Ker(\Phi_i)$ is a maximal ideal of $\bar F$. Similarly for each $+_j\in BT(V_2)$ we get a maximal ideal $Ker(\Phi_j)$ of $\bar F$. By assumption we know that $\bar{F} \cap Aut(V_1)=\bar{F} \cap Aut(V_2)$, and we know:

\[\begin{array}{lllll}
	 \alpha_1+_\cdot...+_\cdot\alpha_n \in Aut (V_1) &\mbox { if and only if }& \alpha_1+_i...+_i\alpha_n \ne 0 \mbox{ for all } i\in I &\mbox{(by Lemma \ref{autodetby0})}\\
	&\mbox { if and only if }& \alpha_1+_\cdot...+_\cdot\alpha_n \notin Ker(\Phi_i).
\end{array}\]
 Thus  $\bar \alpha \in Aut (V_1)$ if and if only $ \bar \alpha\notin \cup _{i\in I}Ker(\Phi_i)$, so $\cup _{i\in I}Ker(\Phi_i)=\cup _{j\in J}Ker(\Phi_j)$.

As each of the additions $+_i$ is different we must have that the maximal ideals $Ker(\Phi_i)$ are distinct and pairwise coprime (as they are maximal ideals). 

Suppose for contradiction that $BT(V_1) \ne BT(V_2)$ so without loss of generality we have  $+_j\in BT(V_2)\setminus BT(V_1)$. So $Ker (\Phi_j)$ is pairwise coprime to $Ker(\Phi_i)$ for each $i$.  We can therefore apply the Chinese remainder theorem for general rings. That is to say that we have a surjection:
\[\Psi: \bar F \rightarrow \oplus_{i\in I} \frac{\bar{F}}{Ker(\Phi_i)} \oplus \frac{\bar{F}}{Ker(\Phi_j)}. \]

As this is a surjection we have an $x\in \bar F$ such that $\Psi(x)=(a_1,...,a_m,0)$ where $m=|I|$ and each $a_i\ne 0$. That is to say  $x\in Ker(\Phi_j)$, but $x\notin Ker(\Phi_i)$ for each $i\in I$, hence $x\in Ker(\Phi_j)\setminus \cup _{i\in I}Ker(\Phi_i)$. 
However, we have that $Ker(\Phi_j)\subseteq \cup _{i\in I}Ker(\Phi_i)$, giving the required contradiction. 

\QED

\begin{rem}\label{0onallbutoneblock}
We can also use the Chinese remainder theorem similarly to find, for a commutative near-vector space, $V$ (with finite block type), 
an $x={\Sigma _{i=1}^{n}}^\cdot \alpha_i$ such that $x\in \cap _{i\in I}Ker(\Phi_i)\setminus Ker(\Phi_j)$. That is to say we have an $x\in \bar F$ which acts as $0$ on all but one block, on that block $x$ acts as a non-zero member of $F$, thus an automorphism,  $\alpha$ say. We may therefore assume it acts as $1$ on the block (just precompose by the inverse of $\alpha$). Multiplying by such an element then becomes similar to projecting onto a single block.
\end{rem}

This is interesting from a model theory point of view, because for every $\alpha \in \bar{F}$ we can add a sentence to our theory stating whether or not $\alpha$ is an automorphism.  So we can fix the block type of any finite block near-vector space model of that theory. That is to say every near-vector space which is a model will be a direct sum of vector spaces over the same fields. In fact there is a more direct way of proving this. The following lemma was pointed out by a referee, and uses the ideas from Proposition \ref{fixedbyauto} to give a precise description of $\bar F$ when $V$ is a finite block near-vector space. We thank the referee for highlighting this to us.

\begin{prop}\label{fieldprod}
Let $(V,F)$ be a finite block type near-vector space, then $\bar F$ is a finite product of fields.
\end{prop}

\textit{Proof:} Let $V=\oplus_{i=1}^n B_i$ be as in Definition \ref{blocktype}, so each $B_i$ is a vector space over some field $F_i= (F, \circ, +_{u_i})$, and $End(B_i)=F_i$. Now, as in Proposition \ref{fixedbyauto} we use the Chinese remainder theorem for arbitrary rings. In particular we can define for each $i$ the map $\Phi_i: \bar F \rightarrow F_i$ which takes $f\in \bar F$ to its restriction in $B_i$. Let $K_i=Ker (\Phi_i)$ then seeing $\bar F$ as a subset of $End(V)$, we get 
\begin{itemize}
	\item $\bar F/K_i \cong F_i;$
	\item  $\cap_{i=1}^n K_i=\{0\}$ (as if $g\in \cap_{i=1}^n K_i$ it acts as $0$ on all the blocks, thus on the whole of $V$).
\end{itemize}
Thus by the Chinese remainder theorem $ \bar F\cong F_1\times...\times F_n$.

\QED

Assuming finite block type also simplifies the situation for ultraproducts. 
By a similar construction to that in Example \ref{fgnotnvs}  we see that the complete theory of an infinite block near-vector space will have models which are not near-vector spaces. 
However, if we start with a $V$ which has a finite number of blocks, then we will not run into the same problem. Below we use Proposition \ref{fixedbyauto} to show that any ultrapower of such a $V$ will be a near-vector space, and thus all models of its complete theory will be near-vector spaces.

\begin{prop}\label{finiteblocknvs}

Suppose $V$ is a near-vector space with finite block type, so  $V=\oplus_{i=1}^m B_i$. Then every ultrapower of $V$ will be a near-vector space.
\end{prop}

\textit{Proof:} In this case we have that:
 \[Q(V)= (B_1\oplus 0\oplus ...\oplus 0) \cup (0\oplus B_2\oplus ...\oplus 0)\cup ...\cup (0\oplus 0\oplus ...\oplus B_m). \]
Let $B_i'=0\oplus ...\oplus B_i\oplus ...\oplus 0$, the near-vector subspace of $V$.

 Let $\mathcal W= \prod V/\mathcal U$ be any (non-principal) ultrapower of $V$. Clearly $\mathcal W$ will be a linear $F$-group, as the axioms are first-order. We then have that:
\[ Q(\mathcal W)= \{w=(w_i)/\mathcal U: \mbox{ almost  all } w_i\in B_j' \mbox{ for some }j\}.\]

Now given an arbitrary element $v=(v_i)/\mathcal U\in \mathcal W$ we have that $v_i=\sum_{j=1}^m b_i^j$ where $b_i^j\in B_j'$. We therefore have that:

\[\begin{array}{lllllll}
	v&=& ((b_1^1,...,b^1_i,....)+...+(b^m_1,...,b^m_i,....))/\mathcal U\\
	&=&\sum _{j=1}^m ((b_i^j)_i/\mathcal U).
\end{array}\]

Now is it clear that each $(b_1^j,...,b^j_i,....)/\mathcal U\in Q(\mathcal W)$ as $b_i^j\in B_j'$ for all $i$. So as each $B_j'$ is a vector space over $F$ we will have a basis in $Q(\mathcal W)$ with the required property.

\QED


\begin{rmk}\label{almostidemp}
Note that, unlike the infinite block case, in Proposition \ref{finiteblocknvs} the quasi-kernel of a finite block near-vector space is definable. To see this first note that we can define each of the blocks. This is because in each block $B_i$ the addition $+_i$ is different from that on the other blocks, in particular for each block $B_i,$ we will have $x_i\in \bar F$, such that $x_i$ acts as zero on all blocks except for $B_i$ (see Remark \ref{0onallbutoneblock}). So we can define the elements of block $B_i$ using the following formula:
\[v\in B_i\mbox{ if and only if }  x_i(v)\ne 0\vee v=0.\]
Now as the quasi-kernel is just the union of the blocks $B_i$, and there are finitely many of these, it is clearly definable (by $\vee_{i=1}^n B_i(v)$). Note that this is quantifier free.

\end{rmk}









As commutative near-vector spaces are modules over $\bar F$, we know from \cite{baur} that we have quantifier elimination to positive primitive formulas (ones of the form $\exists w_1....w_k \bigwedge ^m _{j=1} ( \Sigma  v_i r_{ij} + \Sigma w_l s_{lj}=\bar 0) $ where $r_{ij}, s_{lj} \in \bar F$). In finite block near-vector spaces we can reduce this further to full quantifier elimination. The first step is to prove $\bar F=frac_V(\bar F)$ so that we do not need quantifiers to define inverses of automorphisms.

\begin{lem}\label{F=frafF}
Let $(V,F)$ be a finite block near-vector space, then $\bar F=frac_V(\bar F).$
\end{lem}

\textit{Proof: } 
By Proposition \ref{fieldprod} we have that $\bar F\cong F_1\times...\times F_n$, a product of fields. 
Thus if we take an $f\in Aut (V) \cap \bar F$ it will be of the form $(f_1,...,f_n)$ with each $f_i\in F_i\setminus \{0\}$, so $f_i^{-1}\in F_i$ (as it is a field), thus $f^{-1}=(f_1^{-1},...,f_n^{-1})\in F_1\times...\times F_n \cong \bar F$. So the inverse of any automorphism in $\bar F$ is already an element of $\bar F$.

\QED

Next we show that each block in a finite block near-vector space is quantifier free definable.

\begin{lem}\label{sumsofblockQF}
Let $V=\oplus _{i\in I}B_i$ be a commutative near-vector space over $F$ with $|I|$ finite, where $B_i$ are the blocks of $V$. Then for $\Delta \subseteq I$, the set $\oplus _{i\in \Delta}B_i\subseteq V$ is definable using a quantifier free formula. 
\end{lem}
\textit{Proof:} By Remark \ref{0onallbutoneblock} we have $x_j\in \bar F$ such that $x_j$ acts as $0$ on all $B_i$ where $i\ne j$, and $x_j$ acts as $\gamma_j\ne 0$ on $B_j$. Now $\gamma_j\in F\setminus \{0\}\subseteq Aut(V)$, so $\gamma_j^{-1}$ acts as an inverse to $\gamma_j$ on the whole of $V$. 

Consider $x_j'=x_j\circ \gamma_j^{-1}\in \bar F$, this acts as $0$ on $B_i$ for $i\ne j$ and identity on $B_j$ (think of $x_j'$ as an idempotent). If we consider the sum $\Sigma_{j\in \Delta} x_j'\in \bar F$ then this will act as identity on all blocks $B_j$ for $j\in \Delta$ and $0$ everywhere else. So $v \in \oplus _{i\in \Delta}B_i$ if and only if $\Sigma_{j\in \Delta} x_j'(v)=v$, clearly a quantifier free statement.

\begin{prop}
Let $V$ be a near-vector space over $F$ with finitely many blocks. Then $Th(V)$ has quantifier elimination in the language $\mc L_{Fnvs} $.
\end{prop}

\textit{Proof:} We will use the back and forth method to show this. The steps are as follows:

\begin{enumerate}
	\item There are finitely many blocks in $V$, let's say $V=\oplus _{i=1}^n B_i$, where each $B_i$ is a vector space over $F_i=(F,+_i, \circ)$. Then  the $\omega$-saturated models of $Th(V)$ are $W=\oplus_{i=1}^n W_i$ where each $W_i$ is an infinite-dimensional vector space over $F_i$. The $\omega$-saturated models cannot have finite dimensional blocks as then all vectors in that block would be in the span of finitely many vectors, thus the type stating the vector is not in the span of any finite subset of vectors is not realised.
	
	\item \textit{Back and forth essentially now works as it does in vector spaces.} Let $W=\oplus _{i=1}^n W_i$ and $U=\oplus _{i=1}^n U_i$ be $\omega$-saturated models of $Th(V)$, where $U_i$ and $W_i$ are blocks defined by the same formula (by Remark \ref{almostidemp} this is a quantifier free formula) thus any partial isomorphism from $U$ to $W$ restricts to a partial isomorphism from $U_i$ to $W_i$. 	
	Suppose we have $\bar a = (a^1,...,a^m) \in W^m$, $\bar b =(b^1,...,b^m) \in U^m$ such that $\bar a\equiv \bar b$ (i.e. there is a partial isomorphism between them). Let:
	
	\[\begin{array}{llll}
			a^i=(w_1^i,...,w_n^i) \mbox{ with } w_j^i \in W_j,\\
			b^i=(u_1^i,...,u_n^i) \mbox{ with } u_j^i \in U_j.\\
	\end{array}\]
	
Note that by Lemma \ref{sumsofblockQF} we must have a partial isomorphism  between the blocks $U_j$ and $W_j$, thus from $(w_j^1, ..., w_j^m)$ to $(u_j^1, ..., u_j^m)$. As these blocks are vector spaces over $F_j$, and these have quantifier elimination, $(w_j^1, ..., w_j^m)$ and $(u_j^1, ..., u_j^m)$ must in fact have the same type. 
	
	Consider $c=(c_1,...,c_n)\in W$ with $c_i\in W_i$, we have the following possibilities:
	\begin{enumerate}
		\item $c_j\in span (w_j^1,..., w_j^m)$ in which case we let $d_j$ be the appropriate sum of $u_j^i$'s.
		\item $c_j\not\in span (w_j^1,..., w_j^m)$ in which case we let $d_j \not \in span (u_j^1,..., u_j^m)$.
	\end{enumerate}
	
	We then let $d=(d_1,...,d_n)$, clearly we will have a partial isomorphism between $\bar ac$ and $\bar b d$. 
	
	\item The back part of this back and forth argument works similarly.\QED
		\end{enumerate}
		
\begin{rmk}
One could also prove quantifier elimination by showing that for a finite block near-vector space $(V,F)$ the ring $frac_V(F)$ is von Neumann regular. By \cite{Prest}, modules over commutative	 von Neumann rings have quantifier elimination. We thank Lorna Gregory for pointing this out to us.
 
\end{rmk}

We can use quantifier elimination to show where these theories stand in the model theoretic universe. Below we establish that the complete theory of a near-vector space with finitely many blocks is totally transcendental (Theorem \ref{fnvstt}), and that they have Morley rank the number of blocks, degree $1$. Of course the latter implies the former, but it was thought it nice to give a direct proof of totally transcendentality.

\begin{thm}\label{fnvstt}
Let $V$ be a near-vector space over $F$ with $n$ blocks, then $T=Th(V)$ is totally transcendental.
\end{thm}

\textit{Proof:} Suppose $\mathcal V$ is a saturated model of $T$. Let $|F|=\kappa$, then clearly $|T|=max\{\kappa, \aleph_0\}=\kappa$, so we need to show that over a set of parameters $A$, with $|A|<\kappa$ we have at most $\kappa$ types in one variable. We simply count the possibilities for an element $u\in V$. 

Firstly, we can consider the trace of $u$ in each block. Here, to get the trace we can multiply by an element of $\bar F$ that acts as the identity on that block and $0$ on all others (see Remark \ref{0onallbutoneblock}). This trace will either be:
\begin{enumerate}
	\item in the trace of an element of $span(A)$ - $\kappa$-many possibilities (remember each block is infinite dimensional).
	\item not in the trace of an element of $span(A)$ - one possibility.
\end{enumerate}

Thus, with $n$ blocks, there are $(\kappa+1)^n=\kappa$-many possibilities. 
There are therefore $\kappa$-many one-types over $A$, proving that $T$ is totally transcendental. \QED

\begin{thm}\label{fnvsRM1dMn}
Let $V$ be a near-vector space over $F$ with $n$ blocks, then $Th(V)$ has Morley rank $n$ and Morley degree $1$.
\end{thm}

\textit{Proof:}
Let $V=B_1\oplus ...\oplus B_n$, consider the induced structure on each block, these will be vector spaces over $F_i$ in the standard language, so we have $RM (B_i)=1$ and $dM(B_i)=1$. As $V$ is the direct sum of (finitely many of) these blocks, which are individually definable, the whole structure must have Morley rank equal to that of the blocks (i.e. $RM(V)=n$) and Morley degree $1$.
\QED

The case where $V$ has infinitely many blocks is interesting. Clearly this can only happen when 
$(F, \circ)$ has infinitely many multiplicative automorphisms, each `twist' inducing a new addition on $F$. 
 It is also clear that an infinite block near-vector space $(V, F)$ will still be an $\bar{F}$-module, and thus must be stable and have quantifier elimination to p.p. formulas. However, the case of the finite block near-vector space suggests we could do much better.

\begin{qu}
Given $V$ a near-vector space over $F$ with infinitely many blocks, what are the models of $Th(V)$ in this language? Is the only obstruction to them being near-vector spaces the fact that they would have infinite support? Do we have $QE$ for some expansion of the language? Where does this theory fit into the model theoretic universe? What property does the ring $\bar F$ have in this case? Are they Von Neumann regular? 
\end{qu}

\bibliographystyle{plain}
\bibliography{referencesBHK1}

\end{document}